\documentclass[12pt]{amsart}
\usepackage{amssymb}
\usepackage{amsfonts}
\usepackage{color}

\newcommand{\nc}{\newcommand}
\nc{\thref}[1]{Theorem~\ref{theo:#1}}
\nc{\selabel}[1]{\label{sect:#1}}
\nc{\seref}[1]{Section~\ref{sect:#1}}
\nc{\lelabel}[1]{\label{lemm:#1}}
\nc{\leref}[1]{Lemma~\ref{lemm:#1}}
\nc{\prlabel}[1]{\label{prop:#1}}
\nc{\prref}[1]{Proposition~\ref{prop:#1}}
\nc{\colabel}[1]{\label{coro:#1}}
\nc{\coref}[1]{Corollary~\ref{coro:#1}}
\nc{\exlabel}[1]{\label{exam:#1}}
\nc{\exref}[1]{Example~\ref{exam:#1}}
\nc{\delabel}[1]{\label{defi:#1}}
\nc{\deref}[1]{Definition~\ref{defi:#1}}
\nc{\eqlabel}[1]{\label{equa:#1}}
\nc{\relabel}[1]{\label{rema:#1}}
\nc{\reref}[1]{Lemma~\ref{rema:#1}}
\providecommand{\operatorname}[1]{\mathrm{#1}\,}
\nc{\Hom}{\operatorname{Hom}}
\nc{\Mor}{\operatorname{Mor}}
\nc{\Aut}{\operatorname{Aut}}
\nc{\Ann}{\operatorname{Ann}}
\nc{\Ker}{\operatorname{Ker}}
\nc{\Trace}{\operatorname{Trace}}
\nc{\Char}{\operatorname{Char}}
\nc{\Mod}{\operatorname{Mod}}
\nc{\End}{\operatorname{End}}
\nc{\Spec}{\operatorname{Spec}}
\nc{\Span}{\operatorname{Span}}
\nc{\sgn}{\operatorname{sgn}}
\nc{\Id}{\operatorname{Id}}
\nc{\Com}{\operatorname{Com}}
\nc{\rank}{\operatorname{rank}}

\let\dirsum=\oplus

\let\:=\colon

\newtheorem{de}{Definition}[section]
\newtheorem{lm}[de]{Lemma}
\newtheorem{pr}[de]{Proposition}
\newtheorem{co}[de]{Corollary}
\newtheorem{re}[de]{Remark}
\newtheorem{res}[de]{Remarks}
\newtheorem{te}[de]{Theorem}
\newtheorem{ex}[de]{Example}
\newtheorem{exs}[de]{Examples}

\def\bex{\begin{ex}}
\def\eex{\end{ex}}
\def\bexs{\begin{exs}}
\def\eexs{\end{exs}}
\def\bl{\begin{lm}}
\def\el{\end{lm}}
\def\bc{\begin{co}}
\def\ec{\end{co}}
\def\bt{\begin{te}}
\def\et{\end{te}}
\def\bpr{\begin{pr}}
\def\epr{\end{pr}}
\def\br{\begin{re}}
\def\er{\end{re}}
\def\brs{\begin{res}}
\def\ers{\end{res}}
\def\bd{\begin{de}}
\def\ed{\end{de}}
\def\be{\begin{equation}}
\def\ee{\end{equation}}
\def\bea{\begin{eqnarray*}}
\def\eea{\end{eqnarray*}}
\def\bp{\begin{proof}}
\def\ep{\end{proof}}

\def\l{{\lambda}}
\def\L{{\Lambda}}
\def\g{{\gamma}}
\def\G{{\Gamma}}
\def\d{{\delta}}
\def\D{{\Delta}}

\def\ZZ{{\mathbb Z}}
\def\NN{{\mathbb N}}

\let\dirsum=\oplus

\let\:=\colon

\begin{document}
\title{Unramified Reductors of Filtered and Graded Algebras}


\begin{abstract}
It is a fairly known fact that most of the algebras appearing in the
theory of rings of differential operators, quantized algebras of
different kinds (including many quantum groups), regular algebras in
projective non-commutative geometry, etc... come equipped with a
natural gradation or filtration controlled by some finite
dimensional vector space(s), e.g. the degree one part of filtration
or gradation. In this note we relate the valuations of the algebras
considered to unramified sub-lattices in some vector space(s).

\end{abstract}

\author{Cornel Baetica and Freddy Van Oystaeyen}

\thanks{2000 \textit{Mathematics Subject Classification}. Primary 16W60;
Secondary 16W35.}
\thanks{Research supported by the bilateral project
``New Techniques in Hopf Algebras and Graded Rings" of the Flemish
and Romanian governments. The first author greatfully acknowledges
the support of the University of Antwerp.}

\keywords{Valuation filtration, unramified reduction.}
\date{}
\maketitle

\section*{Introduction}


\noindent In the attempt to finding ``good" reductors of filtered
and/or graded algebras we found out that the ones which satisfy an
{\it unramifiedness} condition are in fact the most suitable for our
purposes, among them the most important fact being that they give
rise to a separated filtration that comes form the given valuation
on the ground field. Starting form this observation we introduce the
notion of the unramified $F$-reductor, respectively unramified
graded reductor and reduce the problem of finding an extending
non-commutative valuation to finding a reductor in an associated
graded ring having a domain for its reduction. As an application we
single out the case of affine algebras, respectively positively
graded algebras although our results apply to a large class of
examples.

Throughout this paper we consider $K$ as being a field, $\G$ a
totally ordered group and $v:K\longrightarrow \G\cup\{\infty\}$ a
valuation on $K$. We always assume that $v$ is surjective, hence
$\G$ is a commutative group. Set $O_v$ the valuation ring of $K$
associated to $v$, $m_v$ its unique maximal ideal and $k_v=O_v/m_v$
the residue field. For a given ring $R$ a family $FR$ of additive
subgroups $F_{\g}R$, $\g\in\G$ satisfying
\begin{itemize}
\item[$(i)$] $\g\leq\d$ implies $F_{\g}R\subseteq F_{\d}R$,
\item[$(ii)$] $F_{\g}R\;F_{\d}R\subseteq F_{\g+\d}R$, for all
$\g,\d\in\G$, \item[$(iii)$] $1\in F_0R$,
\end{itemize}
is a called a {\em $\G$-filtration} on $R$. For a $\G$-filtration
$FR$ we may define the associated graded ring
$G_F(R)=\dirsum_{\g\in\G}F_{\g}R/F_{<\g}R$, where
$F_{<\g}R=\sum_{\g'<\g}F_{\g'}R$. We say that $FR$ is {\em
$\G$-separated} if for every $a\in R$, $a\neq 0$ there is a
$\g\in\G$ such that $a\in F_{\g}R-F_{<\g}R$. Note that for $\G=\ZZ$
the $\G$-separatedness is equivalent to the separatedness, that is
$\cap_{n\in\ZZ}F_nR=0$, but for an arbitrary $\G$ the latter
condition may be strictly weaker. If $FR$ is $\G$-separated, then we
may define the {\em principal symbol map} $\sigma:R\longrightarrow
G_F(R)$ by $\sigma(a)=a$ mod $F_{<\g}R$ whenever $a\in
F_{\g}R-F_{<\g}R$. The {\em degree} $\deg\sigma(a)=\g$ of
$\sigma(a)$ is uniquely determined. To $FR$ one associates a
value-function $v_F:R\longrightarrow \G\cup\{\infty\}$ defined by
$v_F(0)=\infty$ and $v_F(a)=-\deg(a)$ for $a\neq 0$. It is well
known that $v_F$ is a valuation function on $R$ whenever $G_F(R)$ is
a domain; see \cite[Corollary 4.2.7]{VO}. On the other side, for a
valuation $v$ on $K$ as before we consider a $\G$-filtration $f^vK$
on $K$ given by $f^v_{\gamma}K=\{x\in K: v(x)\geq -\gamma\}$.
Obviously $f^v_0K=O_v$ and $f^v_{<0}K=m_v$. All the filtrations
considered in this paper are supposed to be {\it exhaustive}, that
is $\cup_{\g\in\G}F_{\g}R=R$. For other unexplained notations the
reader is referred to \cite{VO}.

\section{Unramified reductions of filtered and graded algebras}

\bd If $\D$ is a skewfield and $\L\subset\D$ is a subring, then $\L$
is a {\em valuation ring} of $\D$ if it is invariant under inner
automorphisms of $\D$ and for $x\in\D$, $x\neq 0$, either $x\in\L$
or $x^{-1}\in\L$.\ed


Let $V$ be a finite dimensional $K$-vector space and $M$ an
$O_v$-submodule of $V$. Recall that $M$ is an {\em $O_v$-lattice} in
$V$ if it contains a $K$-basis of $V$ and is a submodule of a
finitely generated $O_v$-submodule of $V$. Usually, we will denote
the $k_v$-vector space $M/m_vM$ by $\overline{V}$ and the residue
field $k_v$ by $\overline{K}$. For any $O_v$-lattice $M$ in $V$ we
have $\dim_{\overline{K}}\overline{V}\leq\dim_KV$. When equality
holds we say that $M$ defines an {\em unramified reduction} of $V$.

\bpr \label{free lattice}
Let $V$ be a finite dimensional $K$-vector space and $M$
an $O_v$-lattice in $V$.
\begin{itemize}
\item[$(i)$] If $M$ is a finitely generated $O_v$-module,
then $M$ is a free $O_v$-module of rank less or equal than $\dim V$.
(This happens, for instance, when $\G=\ZZ$.)
\item[$(ii)$] If $M$ defines an unramified reduction of $V$,
then $M$ is a free $O_v$-module of rank equal to $\dim V$.
\end{itemize}
\epr \bp (i) It is a well-known fact that finitely generated
torsion-free modules over valuation rings are free. Since the
$O_v$-module $M$ is contained in
a $K$-vector space, it is torsion-free. \\
(ii) Choose a $\overline{K}$-basis
$\{\overline{x}_1,\dots,\overline{x}_n\}$ in $\overline{V}$,
$n=\dim_KV$. The elements $x_1,\dots,x_n\in M$ are linearly
independent over $K$, in particular over $O_v$. Now we show that
$x_1,\dots,x_n$ is a system of generators for $M$ as an
$O_v$-module. (Note that we can not apply Nakayama's Lemma because
we do not know that $M$ is a finitely generated $O_v$-module.) For
an element $x\in M$ there exist $a_1,\dots,a_n\in K$ such that
$x=\sum_{i=1}^na_ix_i$. There also exists an element $a_j\neq 0$
such that $a_i^*=a_j^{-1}a_i\in O_v$ for all $i$. If $a_j\in O_v$,
then all the elements $a_i$ belong to $O_v$, henceforth $x\in
O_vx_1+\cdots+O_vx_n$. Otherwise, $a_j^{-1}\in m_v$ and by taking
the residues modulo $m_vM$ we get
$0=\sum_{i=1}^n\widehat{a}_i^*\overline{x}_i$, a contradiction since
we have $a_j^*=1$. \ep

The {\em valuation filtration} $F^vV$ on $V$ defined by
$F^v_{\g}V=(f^v_{\g}K)M$, $\g\in\G$, is an exhaustive filtration, but {\em we
do not know if it is $\G$-separated or not.} As we see bellow the
valuation filtration is $\G$-separated whenever $M$
is a free $O_v$-module.

\bl \label{separated} If $M$ is a free $O_v$-module, then the
valuation filtration $F^vV$ on $V$ is $\G$-separated. \el \bp Let
$\{x_1,\dots,x_r\}$ be an $O_v$-basis of $M$. Then it is easily seen
that $\{x_1,\dots,x_r\}$ are linearly independent over $K$. Take an
element $x\in V$. From $KM=V$ it follows that there exists an $a\in
O_v-\{0\}$ such that $ax\in M$. Write $ax=\sum_{i=1}^ra_ix_i$ with
$a_i\in O_v$ for all $i$. If we consider $a_i^*=a^{-1}a_i\in K$,
then we can write $x=\sum_{i=1}^ra_i^*x_i$. For those $a_i^*\neq 0$
we denote $v(a_i^*)$ by $-\g_i$, where $\g_i\in\G$ for all $i$.
Assume that $\g_1\leq\dots\leq\g_r=\g$ and show that $x\in
F_{\g}^vV-F_{<\g}^vV$. Obviously $x\in F_{\g}^vV$, and if we suppose
that $x\in F_{<\g}^vV$, then there exists an $\d<\g$ such that
$x\in(f_{\d}^vK)M$. Now we can write $x=\sum_{i=1}^rb_ix_i$ with
$b_i\in f_{\d}^vK$ for all $i$. This entails that $b_i=a_i^*$ for
all $i$, in particular $b_r=a_r^*$, and thus we get $v(b_r)=-\g$. As
$b_r\in f_{\d}^vK$ we have $-\g\geq-\d$, that is $\d\geq\g$, a
contradiction. \ep

Now let us consider $A$ a $K$-algebra and $FA$ an exhaustive
separated $\ZZ$-filtration on $A$ such that $K\subseteq F_0A$.
Suppose in addition that $FA$ is {\em finite}, i.e.
$\dim_KF_nA<\infty$ for all $n\in\ZZ$. Since $FA$ is finite and
separated, it must be {\em left limited}, i.e. there is $n_0\in\ZZ$
such that $F_nA=0$ for all $n\leq n_0$. Without loss of generality
we may suppose that the filtration $FA$ is {\em positive}, that is
$F_nA=0$ for all $n<0$. Consider $\L\subset A$ a subring. The {\em
induced filtration} $F\L$ of $FA$ on $\L$ is given by $F_n\L=\L\cap
F_nA$, $n\in\NN$. Then $\L$ is an {\em $F$-reductor} of $A$ if
$\L\cap K=O_v$ and $F_n\L$ is an $O_v$-lattice in $F_nA$ for all
$n\in\NN$. We call the ring $\overline{A}=\L/m_v\L$ the (filtered)
{\em reduction} of $A$ with respect to $\L$. The {\em valuation
filtration} $F^vA$ on $A$ is defined by $F^v_{\g}A=(f^v_{\g}K)\L$,
$\g\in\G$.

\bd
Let $A$ be a filtered $K$-algebra with a finite filtration
$FA$ and $\L\subset A$ an $F$-reductor. We say that $\L$ is an
{\em unramified $F$-reductor} if $F_n\L$ is an unramified
reduction of $F_nA$ for all $n\in\NN$.
\ed

The existence of (unramified) $F$-reductors in the general case
seems to be unlikely, but in case the algebra is given by a finite
number of generators and finitely many relations that reduce well,
it is easy to find one. Note that finite dimensional algebras over
fields have always unramified reductors; see \cite[Proposition
1.2]{VO1}. When such a reductor there exists the valuation
filtration on $A$ turns out to be $\G$-separated.

\bpr \label{separ. filt.} Let $A$ be a filtered $K$-algebra with a
finite filtration $FA$ and $\L\subset A$ an unramified $F$-reductor.
Then the valuation filtration $F^vA$ on $A$ is $\G$-separated. \epr
\bp $F_n\L$ is a free $O_v$-module for all $n\in\mathbb{N}$. From
Lemma \ref{separated} we get that the valuation filtration
$F^v(F_nA)$ is $\G$-separated for all $n\in\mathbb{N}$, and by using
\cite[Proposition 3.5(ii)]{BVO} we have that $F^vA$ is
$\G$-separated. \ep

Note that for $\G=\ZZ$ an $F$-reductor necessarily defines a
separated filtration. In this case $F_n\L$ is a finitely generated
$O_v$-module, hence it is free; see Proposition \ref{free
lattice}(i).

In the following we prove some properties of the unramified
$F$-reductors. First we show that the unramified $F$-reductors are
free $O_v$-modules.

\bl
Let $\L$ be an $F$-reductor of $A$. Then we have $m_vF_j\L\cap F_i\L=m_vF_i\L$
for all $i,j\in\mathbb{N}$.
\el
\bp
The non-trivial case is $i<j$. In the proof of Proposition 3.5 from \cite{BVO}
we showed that
\begin{equation}\label{1}
(f^v_{\g}K)\L\cap F_nA=(f^v_{\g}K)(\L\cap F_nA)
\end{equation}
for all $\g\in\G$ and $n\in\NN$. From (\ref{1}) we can easily deduce that
$m_v(\L\cap F_nA)=m_v\L\cap F_nA$. So $m_vF_j\L\cap F_i\L=
m_v(\L\cap F_jA)\cap F_i\L=(m_v\L\cap F_jA)\cap F_i\L=m_v\L\cap F_iA\cap\L=
m_v\L\cap F_iA=m_v(\L\cap F_iA)=m_vF_i\L$.
\ep

\bpr Let $A$ be a filtered $K$-algebra with a finite filtration $FA$
and $\L\subset A$ an unramified $F$-reductor. Then $\L$ is a free
$O_v$-module. \epr \bp Let
$\{\overline{x}_1,\dots,\overline{x}_{p_0}\}\subset\overline{F_0A}
=F_0\L/m_vF_0\L$ be a $\overline{K}$-basis. Then
$\{x_1,\dots,x_{p_0}\}$ is an $O_v$-basis for $F_0\L$ and a
$K$-basis for $F_0A$. Since $m_vF_1\L\cap F_0\L=m_vF_0\L$ we get
that $\overline{F_0A}$ is a $\overline{K}$-vector subspace of
$\overline{F_1A}$. Now we extend the $\overline{K}$-basis
$\{\overline{x}_1,\dots, \overline{x}_{p_0}\}$ of $\overline{F_0A}$
to a $\overline{K}$-basis
$\{\overline{x}_1,\dots,\overline{x}_{p_0},\dots,\overline{x}_{p_1}\}$
of $\overline{F_1A}$. Again we have that $\{x_1,\dots,x_{p_0},\dots,
x_{p_1}\}$ is an $O_v$-basis for $F_1\L$ and a $K$-basis for $F_1A$.
We can continue this way and we obtain a sequence $(x_n)_{n\geq 1}$
of elements in $\L$ that forms an $O_v$-basis of $\L$. \ep

The next result shows that the unramifiedness is preserved by taking the
intersection with a filtered subring.

\bpr Let $A$ be a filtered $K$-algebra with a finite filtration
$FA$, $A'$ a $K$-subalgebra of $A$ and $\L\subset A$ an unramified
$F$-reductor. Then $\L'=\L\cap A'$ is an unramified $F$-reductor of
$A'$. \epr \bp We consider on $A'$ the induced filtration
$FA'=FA\cap A'$, and similarly for $\L'$ we have $F\L'=\L'\cap
FA'=\L\cap A'\cap FA=F\L\cap A'$. We have to prove that $F_n\L'$ is
an unramified reduction of $F_nA'$ for all $n\in\mathbb{N}$. It is
enough to show that $\dim_{\overline{K}}
\overline{F_nA'}=\dim_{K}F_nA'$. In order to do that, let us first
note that $m_v(\L\cap F_nA\cap A')=m_v\L\cap F_nA\cap A'$, or
equivalently $m_vF_n\L\cap F_n\L'=m_vF_n\L'$. It follows that
$\overline{F_nA'}$ is a $\overline{K}$-vector subspace of
$\overline{F_nA}$. Now we choose a $\overline{K}$-basis
$\{\overline{x}_1,\dots,\overline{x}_r\}$ in $\overline{F_nA'}$ and
check that the representatives $\{x_1,\dots,x_r\}$ are also a
$K$-basis of $F_nA'$. Let us first extend
$\{\overline{x}_1,\dots,\overline{x}_r\}$ to
$\{\overline{x}_1,\dots,\overline{x}_s\}$, $s\geq r$, a
$\overline{K}$-basis in $\overline{F_nA}$. Since $F_n\L$ is an
unramified reduction of $F_nA$, we have that $\{x_1,\dots,x_s\}$ is
a $K$-basis of $F_nA$. Assume that $\{x_1,\dots,x_r\}$ is not a
$K$-basis of $F_nA'$, and extend it to $\{x_1,\dots,x_r,
y_1,\dots,y_t\}$ a $K$-basis of $F_nA'$. Now we can write
$y_u=\sum_{j=1}^ra_{uj}x_j+\sum_{k=r+1}^sb_{uk}x_k$, where $1\leq
u\leq t$, $a_{uj}\in K$, and $b_{uk}\in K$ not all zero. Set
$z_u=\sum_{k=r+1}^sb_{uk}x_k$, $1\leq u\leq t$. Obviously, $z_u\in
F_nA'$ for all $u=1,\dots,t$, and moreover $\{x_1,\dots,x_r,
z_1,\dots,z_t\}$ is a $K$-basis of $F_nA'$. Let $b\in K$, $b\neq 0$,
such that $b_{uk}^*=bb_{uk}\in O_v$ for all $u=1,\dots,t$,
$k=r+1,\dots,s$, $b_{uk_0}^*=1$ for a $k_0\in\{r+1,\dots,s\}$, and
set $z_u^*=bz_u$. The elements $z_u^*$ belong to $F_nA'\cap
F_n\L=F_n\L'$, and
$\overline{z}_u^*=\sum_{k=r+1}^s\widehat{b}_{uk}^*\overline{x}_k$.
On the other hand, $\overline{z}_u^*\in\overline{F_nA'}$ implies
that $\overline{z}_u^*= \sum_{i=1}^r\widehat{c}_i\overline{x}_i$ for
some $\widehat{c}_i\in \overline{K}$, a contradiction. \ep

The unramifiedness also behaves well with respect to the tensor
products.

\bpr \label{tensor} Let $A$, $A'$ be filtered $K$-algebras with
finite filtrations $FA$, respectively $FA'$, and $\L\subset A$,
respectively $\L'\subset A'$ unramified $F$-reductors. Then
$\L\otimes_{O_v}\L'$ is an unramified $F$-reductor of $A\otimes_KA'$
with respect to the tensor filtration. \epr \bp Note that
$A=K\otimes_{O_v}\L$, and $A'=K\otimes_{O_v}\L'$. This entails that
$A\otimes_{K}A'=(K\otimes_{O_v}\L)\otimes_{K}(K\otimes_{O_v}\L')=
K\otimes_{O_v}(\L\otimes_{O_v}\L')$, and thus we can deduce that
$\L\otimes_{O_v}\L'$ is a subring of $A\otimes_{K}A'$. Now we define
on the tensor product $A\otimes_{K}A'$ a finite filtration (called
{\em tensor} filtration) given by $F_n(A\otimes_{K}A')=
\bigoplus_{i+j=n}F_iA\otimes_KF_jA'$, an similarly
$F_n(\L\otimes_{K}\L')= \bigoplus_{i+j=n}F_i\L\otimes_{O_v}F_j\L'$.
Since $K\otimes_{O_v}F_n\L=F_nA$ and $K\otimes_{O_v}F_n\L'=F_nA'$
for all $n$, we get that $F_i\L\otimes_{O_v}F_j\L'$ is an
$O_v$-submodule of $F_iA\otimes_KF_jA'$ for all $i,j$, and moreover
the filtration defined on $\L\otimes_{O_v}\L'$ is induced by the
filtration defined on $A\otimes_KA'$. The unramifiedness of
$\L\otimes_{O_v}\L'$ is easily seen. \ep

We get now the following well-known result (see \cite[Proposition
2.1]{VO1})

\bc If $A$ is a finitely dimensional $K$-algebra, $A'$ a $K$-central
simple subalgebra of $A$, and $\L'$ an unramified reduction of $A'$,
then there exists $\L$ an unramified reduction of $A$ such that
$\L'=\L\cap A'$. \ec \bp It is a classical result that
$A=A'\otimes_k A''$, where $A''$ is the centralizer of $A'$ in $A$,
and now we apply Proposition \ref{tensor}. \ep

The property of an unramified $F$-reductor of being a valuation ring is completely
described by its reduction over $O_v$.

\bpr Let $\D$ be a skewfield that contains $K$ in its center, $F\D$
a finite filtration on $\D$, and $\L\subset \D$ an unramified
$F$-reductor. Then $\L$ is a valuation ring (for $\D$) if and only
if $\overline{\L}$ is a skewfield. \epr \bp Assume that $\L$ is a
valuation ring for $\D$ with maximal ideal $m$. Obviously
$m_v\L\subset m$, and we aim to show that the foregoing inclusion is
an equality. Pick an element $x\in m$. Then $x^{-1}\in\D-\L$, and
thus we get an $n\in\NN$ such that $x\in F_n\L$ and $x^{-1}\in
F_nA-F_n\L$. Write $x=\sum_{i=1}^ra_ix_i$, $a_i\in O_v$, and
$x^{-1}=\sum_{i=1}^rb_ix_i$, $b_i\in K$ (not all in $O_v$!), where
$\{x_1,\dots,x_r\}$ is an $O_v$-basis for $F_n\L$, and a $K$-basis
for $F_nA$. By standard arguments we get an $j\in\{1,\dots,r\}$ such
that $b_i^*=b_j^{-1}b_i\in O_v$ for all $i$ and $b_j^{-1}\in m_v$.
Then $b_j^{-1}=(\sum_{i=1}^ra_ix_i)(\sum_{i=1}^rb_i^*x_i)$. By
taking the residues modulo $m_vF_n\L$ we get that
$0=(\sum_{i=1}^r\widehat{a}_i\overline{x}_i)
(\sum_{i=1}^r\widehat{b}_i^*\overline{x}_i)$. As $b_j^*=1$, we must
have that $\sum_{i=1}^r\widehat{a}_i\overline{x}_i=0$, and this
implies that $a_i\in m_v$ for all $i$.

Conversely, suppose that $\overline{\L}$ is a skewfield. Since the
valuation filtration $F^v\D$ is $\G$-separated, $F_0^v\D=\L$, and
$G_v(\D)_0=\overline{\L}$ is a domain, we can apply Corollary 1.4
from \cite{BVO} and get that $\L$ is a valuation ring of $\D$. \ep

Let $A$ be an affine $K$-algebra generated by $a_1,\ldots,a_n$,
$K<\underline{X}>$ the free $K$-algebra on the set
$\underline{X}=\{X_1,\ldots,X_n\}$ and
$\pi:K<\underline{X}>\;\longrightarrow A$ the canonical $K$-algebras
morphism given by $\pi(X_i)=a_i$, $i=1,\ldots,n$. Restriction of
$\pi$ to $O_v<\underline{X}>$ defines a subring $\L$ of $A$, i.e.
$\L=\pi(O_v<\underline{X}>)$. As before, the subring $\L$ yields a
valuation filtration $F^vA$ on $A$ given by $F^v_{\g}A=(f^v_{\g}K)\L$,
$\g\in\G$. It is easy to see that for a graded $K$-algebra $A$ that
has a finite $PBW$-basis, the subring $\L$ is an unramified $F$-reductor
with respect to the grading filtration $FA$.

Note now that via Proposition \ref{separ. filt.} we get a new proof
of the following result (Theorem 3.4.7 from \cite{VO})

\bt\label{PBW} For a graded $K$-algebra $A$ that has a finite
$PBW$-basis, the valuation filtration $F^vA$ is $\G$-separated and
strong. \et

By $G_v(A)$ we denote the associated graded ring determined by the
valuation filtration $F^vA$. The next two results can be proved
similarly to their correspondents (Theorem 2.2 and Proposition 2.4)
from \cite{BVO}, but we record them here in order to emphasize that
now we do not need extra-conditions on the filtered parts of $A$.

\bpr\label{extension}
Let $A$ be a $K$-algebra with a finite filtration
and $\L\subset A$ an unramified $F$-reductor of $A$. If
$\overline{A}=\L/m_v\L$ is a domain and $A$ is an
Ore domain, then every valuation $v$ on $K$ extend to
$Q=Q_{cl}(A)$, the classical ring of fractions of $A$.
\epr

\bpr\label{crossed}
If $A$ is a $K$-algebra with a finite filtration
and $\L\subset A$ an unramified $F$-reductor of $A$, then the associated
graded ring $G_v(A)$ is isomorphic to the twisted group ring $\overline{A}*\G$,
where $\overline{A}=\L/m_v\L$.
\epr


It is a common strategy to deduce properties of $A$ from
properties of $G_F(A)$ whenever possible. Let us focus on the
graded situation now.

Let $R$ be an $\NN$-graded $K$-algebra with $K\subseteq R_0$.
Suppose that the gradation is {\em finite}, i.e. $\dim_KR_n<\infty$
for all $n\in\NN$ and let $\L\subset R$ be a graded subring. Then
$\L$ is called a {\em graded reductor} if $\L\cap K=O_v$ and
$\L\cap R_n$ is an $O_v$-lattice in $R_n$ for all $n\in\NN$.
The ring $\overline{R}=\L/m_v\L$ is called the (graded) {\em
reduction} of $R$ with respect to $\L$. The valuation filtration
$F^vR$ on $R$ is similarly defined by $F^v_{\g}R=(f^v_{\g}K)\L$,
$\g\in\G$.


\bd Let $R$ be an $\NN$-graded $K$-algebra with $K\subseteq R_0$, and
$\dim_KR_n<\infty$ for all $n\in\NN$. If $\L\subset R$ is a graded reductor, then
we say that $\L$ is an {\em unramified graded reductor} if $\L\cap R_n$ is an unramified
reduction of $R_n$ for all $n\in\NN$. \ed

It is rather easy to see that we have similar properties for unramified graded
reductors to the ones already proved for unramified $F$-reductors. We mention
here only one of them, but the interested reader is invited to do it on his own.

\bpr \label{separ. gr. filt.} Let $R$ be an $\NN$-graded $K$-algebra
with a finite gradation and $\L\subset R$ an unramified graded
reductor. Then the valuation filtration $F^vR$ on $R$ is
$\G$-separated. \epr \bp $\L\cap R_n$ is a free $O_v$-module for all
$n\in\mathbb{N}$. Lemma \ref{separated} entails that the valuation
filtration $F^vR_n$ is $\G$-separated for all $n\in\mathbb{N}$, and
by using \cite[Proposition 3.7(ii)]{BVO} we have that $F^vR$ is
$\G$-separated. \ep

We also remark that for $\G=\ZZ$ a graded reductor necessarily
defines a separated filtration. We mention another interesting case
when a graded reductor defines a $\G$-separated filtration, the case
of connected positively graded algebras. Recall that a $K$-algebra
$R$ is called a {\em connected} positively graded algebra if
$R=K\dirsum R_1\dirsum\cdots\dirsum R_n\dirsum\cdots=K[R_1]$ and
$\dim_KR_1<\infty$. Let us assume that $\dim R_1=n$,
$\underline{X}=\{X_1,\ldots,X_n\}$ are indeterminates over $K$, and
take $\pi:K<\underline{X}>\;\longrightarrow R$ a presentation of
$R$. If we set $\L=\pi(O_v<\underline{X}>)$, then
$\dim_{\overline{K}}\overline{R}_1=n$ since no elements of degree
one in the gradation of $K<\underline{X}>$ are in $\mathcal{R}$, the
ideal of relations of $R$. Nevertheless, $\dim_KR_n$ and
$\dim_{\overline{K}}\overline{R}_n$ may be different for $n>1$.
However, we can prove the following.

\bpr Let $R$ be connected positively graded $K$-algebra and
$\L\subset R$ defined as before. Then $F^vR$ is $\G$-separated. \epr
\bp Note first that $\L$ is a graded ring, where the gradation is
the one inherited from $O_v<\underline{X}>$ via $\pi$. All we have
to do is to show that $\L$ is a graded subring of $R$, i.e.
$\L_n=\L\cap R_n$ for all $n$. For $n=0$ it means that $O_v=\L\cap
K$ which is obviously true. For $n>0$, pick an element $x\in\L\cap
R_n$. Since $x\in R_n$ there exists an element $a\in O_v$ such that
$ax\in\L_n$, and since $x\in\L$ there exists $y\in
O_v<\underline{X}>$ such that $x=\pi(y)$. Writing $y$ as a sum of
homogeneous components, $y=\sum_{i\geq 0}y_i$, and multiplying the
relation by $a\in O_v$ we get that $a\pi(y_i)=0$ for all $i\neq n$,
therefore $\pi(y_i)=0$ for all $i\neq n$, that is
$x=\pi(y_n)\in\L_n$. Thus we get that $\L$ is a graded reductor and
$\L_n$ is $O_v$-free for all $n$ (note that $\L_n$ is a finitely
generated $O_v$-module), and consequently $F^vR$ is $\G$-separated.
\ep

Note that graded reductor $\L$ defined above is not necessarily
unramified, although all its graded components are free
$O_v$-modules.

For the sake of completeness we recall here the following result
(Lemma 3.3 from \cite{BVO})

\bl\label{lattices}
Let $V$ be a finite dimensional $K$-vector
space, $V'\subset V$ a $K$-vector subspace and $M\subset V$ an
$O_v$-submodule.
\begin{itemize}
\item[$(i)$] If $M$ is an $O_v$-lattice in $V$, then the quotient
module $M/M\cap V'$ is an $O_v$-lattice in $V/V'$.
\item[$(ii)$] If $M/M\cap V'$ is an $O_v$-lattice in $V/V'$ and
$M\cap V'$ is an $O_v$-lattice in $V'$, then $M$ is an
$O_v$-lattice in $V$.
\end{itemize}
\el

\bpr\label{connection}
Let $A$ be a $K$-algebra with a finite filtration $FA$,
and $\L\subset A$ a subring.
\begin{itemize}
\item[$(i)$] If $\L\subset A$ is an unramified $F$-reductor, then $G_F(\L)\subset
G_F(A)$ and $\widetilde{\L}\subset\widetilde{A}$ are unramified graded
reductors.
\item[$(ii)$] If $G_F(\L)\subset G_F(A)$ or
$\widetilde{\L}\subset\widetilde{A}$ are unramified graded reductors, then
$\L\subset A$ is an unramified $F$-reductor.
\end{itemize}
\epr
\bp
(i) In order to show that $G_F(\L)\subset G_F(A)$
is an unramified graded reductor we use Lemma \ref{lattices}(i).
That $\widetilde{\L}\subset\widetilde{A}$ is an
unramified graded reductor follows by the definition of Rees algebra. \\
(ii) By induction using Lemma \ref{lattices}(ii).
\ep

Consequently any unramified $F$-reductor give rise to an unramified graded
reductor. On the other side, an unramified graded reductor is an unramified
$F$-reductor, where $FR$ is the grading filtration.

As applications of the above Proposition \ref{connection} we mention
here two results: the first one is Theorem 2.6 from \cite{PVO}, and
the second one is Proposition 3.2 from \cite{L}. It is worthwhile to
remark that in both cases our results hold for every $\G$-valuation,
while their results are given only for discrete valuations.

\bpr
Let $A$ be a $K$-algebra with a finite filtration $FA$, and $\L\subset A$ a subring
such that $\L\cap K=O_v$ and $K\L=A$.
\begin{itemize}
\item[$(i)$] If $G_F(\L)_n$ is a finitely generated $O_v$-module for all $n$,
then the valuation filtration $F^vA$ is $\G$-separated.
\item[$(ii)$] If $F_n\L$ are finitely generated $O_v$-modules for all $n$, then
the filtrations $F^vA$ and $f^vG_F(A)$ are $\G$-separated.
\end{itemize}
\epr \bp (i) Since $G_F(\L)_n$ are all finitely generated
$O_v$-modules, we easily get that $F_n\L$ are all finitely generated
$O_v$-modules, and thus $\L$ is an $F$-reductor. Once again we apply
Proposition 1.1(i) and deduce that $F_n\L$ are $O_v$-free, and this
is enough in order to conclude that the valuation filtration
$F^vA$ is $\G$-separated. \\
(ii) From (i) we know that $F^vA$ is $\G$-separated. Moreover,
$G_F(\L)\subset G_F(A)$ is a graded reductor with the graded
components free $O_v$-modules, therefore the valuation filtration
$f^vG_F(A)$ is $\G$-separated. \ep

\bpr Let $A=K[a_1,\dots,a_r]$ be an affine $K$-algebra, and $\L$ as
before but with the generator filtration, i.e. $F_{-1}\L=0$,
$F_0\L=O_v$, and $F_n\L$ is the $O_v$-module generated by the
elements $a_{i_1}\cdots a_{i_t}$ with $1\leq t\leq n$. If the
associated graded ring $G_F(\L)$ is a flat $O_v$-module, then the
valuation filtration $F^vA$ is $\G$-separated. \epr \bp We consider
on $A$ the generator filtration $FA$. In principle, we do not know
that $G_F(\L)$ is contained in $G_F(A)$, and that is why we are
asking for flatness. To enter the details, take $x\in F_n\L\cap
F_{n-1}A$. We want to show that $x\in F_{n-1}\L$. Suppose that this
is not true. Then there exists an $a\in O_v-\{0\}$ such that $ax\in
F_{n-1}\L$ which means that $ax=0$ in the associated graded ring
$G_F(\L)$. As we know that $G_F(\L)$ is a flat $O_v$-module, i.e.
torsion-free, we get a contradiction. Therefore $F_n\L\cap
F_{n-1}A=F_{n-1}\L$, that is $G_F(\L)\subset G_F(A)$ is a graded
reductor with finitely generated (hence free) graded components, and
this is enough to see that $\L\subset A$ is an $F$-reductor with all
filtered parts $O_v$-free. \ep

Here there are some (old) examples of algebras that admit unramified
reductors.


\bexs {\rm ($i$) Consider a field $k$ with $\Char k \neq 2$, and set
$K=k(X)$ the field of rational functions over $k$. Let $L=K(\xi)$,
where $\xi$ is a root of the irreducible polynomial $T^2+(X-1)T+X\in
K[T]$. Now consider the valuation ring of $K$ given by
$O_v=k[X]_{(X)}$, and $\L=O_v[\xi]$. It is easy to see that
$O_v\subset\L$ is an integral extension, $\overline{K}=k_v\simeq k$,
$\overline{\L}=\L/m_v\L\simeq k[u]$, where $u$ is an idempotent, and
therefore $\L$ is an unramified reductor of $L$.

\noindent ($ii$) Consider $g$ a finite dimensional {\it Lie algebra}
over a field $K$ and $A=U(g)$ the enveloping algebra of $g$. Let
$O_v$ be a valuation ring of $K$. We define a finite dimensional Lie
algebra $g_{O_v}$ over $O_v$ with the same basis and the induced
bracket. Let us fix a $K$-basis $\{x_1,\ldots,x_n\}$ for $g$. We
have structure constants $\l_{ij}^k\in K$ with
$[x_i,x_j]=\sum_{k=1}^n\l_{ij}^kx_k$. Without loss of generality we
may assume that $\l_{ij}^k\in O_v$ (up to multiplying all $x_i$ by a
suitable constant in $O_v$) but not all in $m_v$. Set
$g_{O_v}=O_vx_1+\cdots+O_vx_n$. This is a Lie $O_v$-algebra with the
induced bracket. Furthermore, $g_{O_v}$ is an $O_v$-lattice in $g$.
On $\overline{g}=g_{O_v}/m_vg_{O_v}$ we define a Lie algebra
structure over $\overline{K}$ by setting
$[\overline{x}_i,\overline{x}_j]=\sum_{k=1}^n\widehat{\l_{ij}^k}
\overline{x}_k$, where the $\overline{x}_i$ are the images of the
$x_i$ in $\overline{g}$ and $\widehat{\l_{ij}^k}$ are the images of
$\l_{ij}^k$ in $\overline{K}$. By our assumptions $\overline{g}$ is
not the trivial Lie algebra. Of course, $\overline{g}$ depends on
the choice of the $K$-basis in $g$.

Let $\L=U_{O_v}(g_{O_v})$ be the enveloping algebra of $g_{O_v}$.
Consider on $A$ the standard filtration $FA$ and on $\L$ the induced
filtration $F\L$. We have that the filtration $FA$ is finite,
$G_F(\L)=O_v[X_1,\ldots,X_n]$ is a subring of the polynomial ring
$G_F(A)=K[X_1,\ldots,X_n]$ and obviously it is an unramified graded
reductor. From Proposition \ref{connection}(ii) we get that $\L$ is
an unramified $F$-reductor. Furthermore, the filtration $F^vA$ is
$\G$-separated, $G_v(A)$ is a domain isomorphic to
$U_{k_v}(\overline{g})$ and thus we can extend $v$ to
$D(g)=Q_{cl}(U(g))$.

\noindent ($iii$) For the {\it Weyl algebra} $R=A_n(K)$ we take
$\L=A_n(O_v)$. We claim that $\L$ is an unramified graded reductor
of $R$. First note that $\L$ is a free $O_v$-module, therefore it
defines a good reduction of $R$. On the other side, we have that
$\rank_{O_v}\L_n=\dim_KR_n$ for all $n\in\NN$, since $\L$ and $R$
have the same $PBW$-basis, and so $\L$ is an unramified graded
reductor of $R$. It is also an $F$-reductor since the associated
graded ring of $A_n(K)$ with respect to the Bernstein filtration is
a polynomial ring.

\noindent ($iv$) Let $K$ and $O_v$ as before. Set
$R=K<X,Y>/(XY-qYX)$ for the {\it quantum plane}, where $q$ is a unit
in $O_v$. Then the subring $\L=O_v<X,Y>/(XY-qYX)$ is an unramified
$F$-reductor with respect to the generator (grading) filtration, and
an unramified graded reductor with respect to the mixed gradation.

\noindent ($v$) Let $K$, $O_v$ and $q$ as before, and set
$A=K<X,Y>/(XY-qYX-1)$ for the {\it quantum Weyl algebra}. Then
$\L=O_v<X,Y>/(XY-qYX-1)$ is an unramified $F$-reductor with respect
to the generator filtration.} \eexs


\vspace*{3mm}
\begin{flushright}
\begin{minipage}{148mm}\sc\footnotesize
University of Bucharest, Faculty of Mathematics, Str.
Academiei 14, RO--010014, Bucharest, Romania\\
{\it E--mail address}: {\tt baetica@al.math.unibuc.ro}
\vspace*{3mm}

University of Antwerp, Department of Mathematics and Computer Science,
Middelheimlaan 1, B--2020 Antwerp, Belgium\\
{\it E--mail address}: {\tt fred.vanoystaeyen@ua.ac.be}
\end{minipage}
\end{flushright}%

\end{document}